\documentclass[11pt]{article}

\newcommand{\qed}{{\hfill\rule{4pt}{7pt}}}

\usepackage{amsmath, epsfig, cite}
\usepackage{amssymb}
\usepackage{amsfonts, color}
\usepackage{latexsym}
\usepackage{amsthm}

\newtheorem{thm}{Theorem}[section]

\newtheorem{cor}[thm]{Corollary}

\newtheorem{lem}[thm]{Lemma}
\newtheorem{conj}[thm]{Conjecture}

\theoremstyle{remark}

\newcommand{\pf}{\noindent{\it Proof.} }

\def\N{{\mathbb N}}
\def\Z{{\mathbb Z}^+}

\numberwithin{equation}{section}

\begin{document}

\begin{center}
{\Large\bf Factors of sums and alternating sums involving\\[5pt] binomial coefficients and powers of integers}
\end{center}

\vskip 2mm \centerline{Victor J. W. Guo$^1$ and Jiang Zeng$^{2}$}
\begin{center}
{\footnotesize $^1$Department of Mathematics, East China Normal University,\\ Shanghai 200062,
 People's Republic of China\\
{\tt jwguo@math.ecnu.edu.cn,\quad http://math.ecnu.edu.cn/\textasciitilde{jwguo}}\\[10pt]
$^2$Universit\'e de Lyon; Universit\'e Lyon 1; Institut Camille
Jordan, UMR 5208 du CNRS;\\ 43, boulevard du 11 novembre 1918,
F-69622 Villeurbanne Cedex, France\\
{\tt zeng@math.univ-lyon1.fr,\quad
http://math.univ-lyon1.fr/\textasciitilde{zeng}} }
\end{center}

\vskip 0.7cm {\small \noindent{\bf Abstract.}
We study divisibility properties of certain sums and alternating sums involving binomial coefficients
and powers of integers. For example, we prove that for all positive integers
$n_1, \ldots, n_m$, $n_{m+1}=n_1$, and any nonnegative integer $r$, there holds
\begin{align*}
\sum_{k=0}^{n_1}\varepsilon^k (2k+1)^{2r+1}\prod_{i=1}^{m} {n_i+n_{i+1}+1\choose n_i-k} \equiv 0 \mod (n_1+n_m+1){n_1+n_m\choose n_1},
\end{align*}
and conjecture that for any nonnegative integer $r$ and positive integer $s$ such that $r+s$ is odd,
$$
\sum_{k=0}^{n}\varepsilon ^k (2k+1)^{r}\left({2n\choose n-k}-{2n\choose n-k-1}\right)^{s}
\equiv 0 \mod{ {2n\choose n}},
$$
where $\varepsilon=\pm 1$.



\section{Introduction}
There have been lasting interests in binomial sums.
Although some binomial sums have no closed formulas, it is still possible to show that  they have  some nice factors.
For example,  a result of Calkin~\cite{Calkin} reads
\begin{align*}
\sum_{k=-n}^n (-1)^k{2n\choose n+k}^m\equiv 0\mod{{2n\choose n}}\qquad \textrm{for $m\geq 1$.}
\end{align*}
Generalizing Calkin's result,
Guo, Jouhet, and Zeng \cite{GJZ} proved, among other things, that
$$
\sum_{k=-n_1}^{n_1}(-1)^k\prod_{i=1}^m
{n_i+n_{i+1}\choose n_i+k}\equiv 0\mod{ {n_1+n_m\choose n_1} }
$$
for all $n_1,\ldots,n_{m}\geq 1$ and $n_{m+1}=n_1$. Recently,
motivated by  the moments of the Catalan triangle (see \cite{CC,GHMR}),
Guo and Zeng~\cite{GZ10} were  led to  study some different  binomial sums and proved the congruence
\begin{align*}
2\sum_{k=1}^{n_1}k^{2r+1}\prod_{i=1}^{m} {n_i+n_{i+1}\choose n_i+k}\equiv 0\mod{ n_1{n_1+n_m\choose n_1} }
\end{align*}
for all $n_1,\ldots,n_{m}\geq 1$ and $n_{m+1}=n_1$.

In this paper we will prove some divisibility properties of another kind of sums and alternating sums
involving binomial coefficients and powers of integers. Let $\N$ denote the set of nonnegative integers and
$\Z$ the set of positive integers.  One of our main results may be stated as follows.

\begin{thm}\label{thm:oddmain2}
For all $n_1,\ldots,n_m\in\Z$, $n_{m+1}=n_1$, and $r\in\N$, there holds
\begin{align*}
\sum_{k=0}^{n_1}\varepsilon^k (2k+1)^{2r+1}\prod_{i=1}^{m} {n_i+n_{i+1}+1\choose n_i-k} \equiv 0 \mod (n_1+n_m+1){n_1+n_m\choose n_1},
\end{align*}
where $\varepsilon=\pm 1$.
\end{thm}
Actually,   we shall  derive Theorem \ref{thm:oddmain2}
from the following more general result.

\begin{thm}\label{thm:oddmain}
For all $n_1,\ldots,n_m\in\Z$, $n_{m+1}=n_1$, and $r\in\N$, there hold
\begin{align*}
&\hskip -3mm \sum_{k=0}^{n_1}k^r(k+1)^r(2k+1)\prod_{i=1}^{m} {n_i+n_{i+1}+1\choose n_i-k} \\
&\equiv 0 \mod (n_1+n_m+1){n_1+n_m\choose n_1}n_1^{\min\{1,r\}}n_m^{\min\{1,{r\choose 2}\}}, \\
&\hskip -3mm \sum_{k=0}^{n_1}(-1)^k k^r(k+1)^r(2k+1)\prod_{i=1}^{m} {n_i+n_{i+1}+1\choose n_i-k}\\
&\equiv 0 \mod (n_1+n_m+1){n_1+n_m\choose n_1}n_1 ^{\min\{1,r\}}n_m ^{\min\{1,r\}}.
\end{align*}
\end{thm}

Indeed,
 by the expansion
$$
(2k+1)^{2r}
=(4k^2+4k+1)^r =\sum_{i=0}^r{r\choose i}4^{i}k^i(k+1)^i,
$$
it is clear that  Theorem \ref{thm:oddmain} infers Theorem \ref{thm:oddmain2}.

Recently, Miana and Romero \cite[Theorem 10]{MR2}  evaluated  the moments
$\Psi_r(n):=\sum_{k=0}^n(2k+1)^r A_{n,k}^2$, where
$A_{n,k}$ ($0\leq k\leq n)$ are the ballot numbers defined by
$$
A_{n,k}=\frac{2k+1}{2n+1}{2n+1\choose n-k}={2n\choose n-k}-{2n\choose n-k-1}.
$$
By  the formulas of $\Psi_{1}(n)$, $\Psi_{3}(n)$ and $\Psi_{5}(n)$ (see \cite[Remark 11]{MR2}),
we conjecture that $\Psi_{2r+1}(n)$ is divisible by ${2n\choose n}$.
More generally, we have
\begin{conj}\label{conj:first}
For all $r\in\N$ and $n,s\in\Z$
such that $r+s\equiv 1\pmod 2$, there holds
\begin{align}
\sum_{k=0}^{n}\varepsilon ^k (2k+1)^{r}A_{n,k}^{s} \equiv 0 \mod{ {2n\choose n}}, \label{eq:ank-2}
\end{align}
where $\varepsilon=\pm 1$.
\end{conj}
As a check,  let  $n=7$,   for  all $r\geq0$ and $s\geq 1$ such that $r+s\equiv 1\pmod 2$,  the sum
\begin{align*}
\sum_{k=0}^{7}(2k+1)^{r}A_{7,k}^{s}=429^s+ 3^r 1001^s+5^r 1001^s+ 7^r 637^s+9^r 273^s+11^r 77^s+13^{r+s}+15^r
\end{align*}
is obviously divisible by
${14\choose 7}=8\times3\times11\times13$.

In this paper, we  shall confirm    Conjecture \ref{conj:first} in some  special cases.
\begin{thm}\label{thm:oddpower}
The congruence \eqref{eq:ank-2} holds if $2n+1$ is a prime power or $s=1$.
\end{thm}

In  the next three sections  we shall  provide the proof of Theorem \ref{thm:oddmain} corresponding, respectively,
to the cases $m=1$, $m=2$ and $m\geq 3$.
We then prove Theorem \ref{thm:oddpower} in Section~5.
Finally we give some further consequences of Theorem~\ref{thm:oddmain2} and related conjectures in Section~6.


\section{Proof of Theorem \ref{thm:oddmain} for $m=1$}
Let
\begin{align*}
P_{r}(n)&:=\sum_{k=0}^n {2n+1\choose n-k}k^r(k+1)^r(2k+1), \\ 
Q_{r}(n)&:=\sum_{k=0}^n (-1)^k {2n+1\choose n-k}k^r(k+1)^r(2k+1).
\end{align*}

The $m=1$ case of Theorem \ref{thm:oddmain} may be stated as follows.
\begin{thm}\label{thm:oddone}
For all $n\in\Z$ and $r\in\N$, there hold
\begin{align*}
&P_{r}(n) \equiv 0 \mod (2n+1){2n\choose n}n^{\min\{2,r\}}, \\
&Q_{r}(n) \equiv 0 \mod (2n+1){2n\choose n}n^{\min\{2,2r\}}.
\end{align*}
\end{thm}
\noindent{\it Proof of Theorem {\rm\ref{thm:oddone}.}}
We proceed by induction on $r$.  For $r=0$, we have
\begin{align*}
P_{0}(n)
&=(2n+1)\sum_{k=0}^n \left({2n\choose n-k}-{2n\choose n-k-1}\right)
=(2n+1){2n\choose n},  \\
Q_{0}(n)
&=(2n+1)\sum_{k=0}^n (-1)^k\left({2n\choose n-k}-{2n\choose n-k-1}\right)
=\begin{cases}
0, &\text{if $n>0$,}\\
1, &\text{if $n=0$.}
\end{cases}
\end{align*}
For $r\geq 1$, observing that
$$
{2n+1\choose n-k}k(k+1)=n(n+1){2n+1\choose n-k}-2n(2n+1){2n-1\choose n-k-1},
$$
we have
\begin{align}
P_{r}(n)&=n(n+1)P_{r-1}(n)-2n(2n+1)P_{r-1}(n-1),  \label{eq:rec1}\\
Q_{r}(n)&=n(n+1)Q_{r-1}(n)-2n(2n+1)Q_{r-1}(n-1)   \label{eq:rec2}
\end{align}
for $n\geq 1$. For the above recurrences we derive immediately that
\begin{align*}
P_{1}(n)  &=n(2n+1){2n\choose n}, \qquad
P_{2}(n)  =2n^2(2n+1){2n\choose n},
\end{align*}
and $Q_{1}(1)=-6$,  $Q_{1}(2)=0$, $Q_{2}(1)=-12$, $Q_{2}(2)=120$, and
$Q_{1}(n) =Q_{2}(n)=0$ for $n\geq 3$.
Therefore, Theorem \ref{thm:oddone} is true for $r=0,1,2$. Now suppose that $r\geq 3$ and Theorem \ref{thm:oddone} holds for
$r-1$. Then  $P_{r-1}(n)$ is divisible by
$(2n+1){2n\choose n}n^{2}$, and $2n(2n+1)P_{r-1}(n-1)$ is divisible by
$$2n(2n+1)(2n-1){2n-2\choose n-1}=(2n+1){2n\choose n}n^2.$$
By \eqref{eq:rec1} we see that $P_{r}(n)$ is also divisible by $(2n+1){2n\choose n}n^2$.  This completes the inductive step
for $P_r(n)$. Similarly, by \eqref{eq:rec1} we can prove the case for $Q_r(n)$. \qed

\medskip
We may also consider the following sums:
\begin{align*}
U_{r}(n)&:=\sum_{k=0}^n {2n+1\choose n-k}(2k+1)^{2r}, \\ 
V_{r}(n)&:=\sum_{k=0}^n (-1)^k {2n+1\choose n-k}(2k+1)^{2r}.
\end{align*}
It is easy to see that
\begin{align}
\sum_{k=0}^n{2n+1\choose n-k}
&=\frac{1}{2}\left(\sum_{k=0}^n{2n+1\choose k}
+\sum_{k=0}^n{2n+1\choose n+1+k} \right)=\frac{1}{2}(1+1)^{2n+1}=4^n,\label{eq:bino1}  \\
\sum_{k=0}^n (-1)^k{2n+1\choose n-k}
&=\sum_{k=0}^n (-1)^k\left(  {2n\choose n-k}+{2n\choose n-k-1}\right)={2n\choose n}.  \label{eq:bino2}
\end{align}
Similarly to \eqref{eq:rec1} and \eqref{eq:rec2}, we have
\begin{align}
U_{r}(n)&=(2n+1)^2 U_{r-1}(n)-8n(2n+1)U_{r-1}(n-1),  \label{eq:rec3}\\
V_{r}(n)&=(2n+1)^2 V_{r-1}(n)-8n(2n+1)V_{r-1}(n-1).  \label{eq:rec4}
\end{align}
By \eqref{eq:bino1}--\eqref{eq:rec4} we immediately obtain the following result.
\begin{cor}For $n\in\Z$ and $r\in\N$, there hold
\begin{align*}
4^{-n}\sum_{k=0}^n {2n+1\choose n-k}(2k+1)^{2r} &\equiv 1 \pmod 2,\\
2^{-\alpha(n)}\sum_{k=0}^n (-1)^k {2n+1\choose n-k}(2k+1)^{2r} &\equiv 1 \pmod 2,
\end{align*}
where $\alpha(n)$ denotes the number of $1$'s in the binary expansion of $n$.
\end{cor}

\section{Proof of Theorem \ref{thm:oddmain} for $m=2$ }

We first give two combinatorial identities.
\begin{lem}\label{lem:n1n2} For all $n_1,n_2\in\N$, there hold
\begin{align}
\sum_{k=0}^{n_1}{n_1+n_2+1\choose n_1-k}{n_1+n_2+1\choose n_2-k}(2k+1)
&=(n_1+n_2+1){n_1+n_2\choose n_1}^2,   \label{eq:n1n2-one}\\
\sum_{k=0}^{n_1}(-1)^k {n_1+n_2+1\choose n_1-k}{n_1+n_2+1\choose n_2-k}(2k+1)
&=(n_1+n_2+1){n_1+n_2\choose n_1}.  \label{eq:n1n2-two}
\end{align}
\end{lem}

\pf It is easy to see that the left-hand side of \eqref{eq:n1n2-one} may be written as
\begin{align*}
&\hskip -3mm
\sum_{k=0}^{n_1}(n_2+k+1){n_1+n_2+1\choose n_1-k}{n_1+n_2\choose n_2-k}  \\
&\qquad{}-\sum_{k=0}^{n_1}(n_2+k+2){n_1+n_2+1\choose n_1-k-1}{n_1+n_2\choose n_2-k-1} \\
&=(n_2+1){n_1+n_2+1\choose n_1}{n_1+n_2\choose n_2},
\end{align*}
which is equal to the right-hand side of \eqref{eq:n1n2-one}.

Replacing $k$ by $-k-1$ in the left-hand side of \eqref{eq:n1n2-two}, we observe that
\begin{align*}
&\hskip -3mm
\sum_{k=0}^{n_1}(-1)^k {n_1+n_2+1\choose n_1-k}{n_1+n_2+1\choose n_2-k}(2k+1) \\
&=\sum_{k=-n_1-1}^{-1}(-1)^k {n_1+n_2+1\choose n_1-k}{n_1+n_2+1\choose n_2-k}(2k+1).
\end{align*}
It follows that \eqref{eq:n1n2-two} is equivalent to
\begin{align}
\sum_{k=-n_1-1}^{n_1}(-1)^k {n_1+n_2+1\choose n_1-k}{n_1+n_2+1\choose n_2-k}(2k+1)
=2(n_1+n_2+1){n_1+n_2\choose n_1}.  \label{eq:n1n2-two2}
\end{align}
Since
\begin{align*}
{n_1+n_2+1\choose n_1-k}{n_1+n_2+1\choose n_2-k}(2k+1)
&=(n_1+n_2+1){n_1+n_2\choose n_1-k}{n_1+n_2+1\choose n_2-k} \\
&\qquad{}-(n_1+n_2+1){n_1+n_2+1\choose n_1-k}{n_1+n_2\choose n_2-k-1},
\end{align*}
to prove \eqref{eq:n1n2-two2}, it suffices to establish the following two identities:
\begin{align*}
\sum_{k=-n_1-1}^{n_1}(-1)^k{n_1+n_2\choose n_1-k}{n_1+n_2+1\choose n_2-k}
&={n_1+n_2\choose n_1}, \\
\sum_{k=-n_1-1}^{n_1}(-1)^{k+1}{n_1+n_2+1\choose n_1-k}{n_1+n_2\choose n_2-k-1}
&={n_1+n_2\choose n_1},
\end{align*}
which follow immediately by comparing the coefficients of $x^{2n_2}$ in the expansion of
\begin{align*}
(1-x)^{n_1+n_2}(1+x)^{n_1+n_2+1}=(1+x)(1-x^2)^{n_1+n_2}.
\end{align*}
This completes the proof.  \qed

\medskip
\noindent{\it Remark.}
 We can give another proof of
 \eqref{eq:n1n2-one} and \eqref{eq:n1n2-two} by computing their  generating functions and using
  the following identity
\begin{align}
\sum_{m,n=0}^{\infty} {m+n+\alpha\choose m}{m+n+\beta\choose n}x^m y^n
=\frac{2^{\alpha+\beta}}{\Delta(1-x+y+\Delta)^\alpha(1+x-y+\Delta)^\beta},\label{eq:delta}
\end{align}
where $\Delta:=\sqrt{1-2x-2y-2xy+x^2+y^2}$.
The identity \eqref{eq:delta} is equivalent to the generating function of Jacobi polynomials.
See \cite[p.~298]{AAR}, \cite[p.~271]{Ra} or \cite{FL} for a proof of this identity, and \cite{GZ06} for an application
to prove  some double-sum binomial coefficient identities.

As an example, we compute the generating function for the left-hand side of \eqref{eq:n1n2-two}:
\begin{align*}
&\hskip -3mm
\sum_{n_1,n_2=0}^\infty x^{n_1}y^{n_2}\sum_{k=0}^{\min\{n_1,n_2\}}(-1)^k {n_1+n_2+1\choose n_1-k}{n_1+n_2+1\choose n_2-k} (2k+1)\\
&=\sum_{k=0}^\infty (2k+1)(-xy)^k \sum_{n_1,n_2=k}^\infty x^{n_1-k}y^{n_2-k} {n_1+n_2+1\choose n_1-k}{n_1+n_2+1\choose n_2-k} \\
&=\sum_{k=0}^\infty \frac{2^{4k+2}(2k+1)(-xy)^k }{\Delta(1-x+y+\Delta)^{2k+1}(1+x-y+\Delta)^{2k+1}}\\
&=\sum_{k=0}^\infty \frac{2^{2k+1}(2k+1)(-xy)^k }{\Delta(1-x-y+\Delta)^{2k+1}}\\
&=\frac{2(1-4xy(1-x-y+\Delta)^{-2})}
{\Delta(1-x-y+\Delta)\left(1+4xy(1-x-y+\Delta)^{-2}\right)^2} \\
&=\frac{1}{(1-x-y)^2}.
\end{align*}
Clearly the last expression  is the generating function for the right-hand side of \eqref{eq:n1n2-two}.

\medskip

Let
\begin{align*}
P_{r}(n_1,n_2)&:=\sum_{k=0}^{n_1}{n_1+n_2+1\choose n_1-k}{n_1+n_2+1\choose n_2-k}k^r(k+1)^r (2k+1),\\
Q_{r}(n_1,n_2)&:=\sum_{k=0}^{n_1}(-1)^k {n_1+n_2+1\choose n_1-k}{n_1+n_2+1\choose n_2-k}k^r(k+1)^r (2k+1).
\end{align*}

Then,  the $m=2$ case of Theorem \ref{thm:oddmain} may be stated as follows.

\begin{thm}\label{thm:oddtwo}
For all $n_1,n_2\in\Z$ and $r\in\N$, there hold
\begin{align*}
P_{r}(n_1,n_2)
&\equiv 0 \mod (n_1+n_2+1){n_1+n_2\choose n_1}n_1^{\min\{1,r\}}n_2^{\min\{1,{r\choose 2}\}}, \\
Q_{r}(n_1,n_2)
&\equiv 0 \mod (n_1+n_2+1){n_1+n_2\choose n_1}n_1 ^{\min\{1,r\}}n_2 ^{\min\{1,r\}}.
\end{align*}
\end{thm}

\noindent{\it Proof of Theorem {\rm\ref{thm:oddtwo}.}}
We proceed by induction on $r$.
Writing
\begin{align*}
{n_1+n_2+1\choose n_1-k}{n_1+n_2+1\choose n_2-k}k(k+1)
&=n_1(n_1+1){n_1+n_2+1\choose n_1-k}{n_1+n_2+1\choose n_2-k}  \\
&\quad{}-(n_1+n_2+1)^2{n_1+n_2\choose n_1-k-1}{n_1+n_2\choose n_2-k},
\end{align*}
we derive
\begin{align}
P_{r}(n_1,n_2)
&=n_1(n_1+1) P_{r-1}(n_1,n_2)-(n_1+n_2+1)^2 P_{r-1}(n_1-1,n_2), \label{eq:recn1n2-1}\\
Q_{r}(n_1,n_2)
&=n_1(n_1+1) Q_{r-1}(n_1,n_2)-(n_1+n_2+1)^2 Q_{r-1}(n_1-1,n_2),\quad r\geq 1.  \label{eq:recn1n2-2}
\end{align}
From the above recurrences and Lemma \ref{lem:n1n2} we immediately get
\begin{align*}
P_{1}(n_1,n_2)&=n_1(n_1+n_2+1){n_1+n_2\choose n_1}{n_1+n_2-1\choose n_1},  \\
P_{2}(n_1,n_2)&=2n_1n_2(n_1+n_2+1){n_1+n_2\choose n_1}{n_1+n_2-2\choose n_1-1},\\
Q_{1}(n_1,n_2)&={}-n_1n_2(n_1+n_2+1){n_1+n_2\choose n_1},  \\
Q_{2}(n_1,n_2)&=n_1n_2(n_1n_2-n_1-n_1-1)(n_1+n_2+1){n_1+n_2\choose n_1}.
\end{align*}
Therefore, Theorem \ref{thm:oddtwo} is true for $r=0,1,2$.
Now suppose that the statement is true for $r-1$ ($r\geq 3$). Then $n_1 P_{r-1}(n_1,n_2)$ is divisible by
$n_1 n_2(n_1+n_2+1){n_1+n_2\choose n_1}$ and $(n_1+n_2+1) P_{r-1}(n_1-1,n_2)$ is divisible by
$$
(n_1+n_2+1) n_2(n_1+n_2){n_1+n_2-1\choose n_1-1}=n_1 n_2(n_1+n_2+1) {n_1+n_2\choose n_1}.
$$
By \eqref{eq:recn1n2-1}, we see that $P_{r}(n_1,n_2)$ is divisible by $n_1 n_2(n_1+n_2+1){n_1+n_2\choose n_1}$, and
we complete the inductive step for $P_{r}(n_1,n_2)$. The case for $Q_{r}(n_1,n_2)$ is exactly the same.
\qed

\medskip
\noindent{\it Remark.}
Similarly, if we set
\begin{align*}
\overline{P}_{r}(n_1,n_2):=\sum_{k=0}^{n_1}{n_1+n_2\choose n_1-k}{n_1+n_2\choose n_2-k}k^{2r+1},
\end{align*}
then
\begin{align*}
\overline{P}_{0}(n_1,n_2)&=\frac{n_1}{2}{n_1+n_2\choose n_1}{n_1+n_2-1\choose n_1}, \\
\overline{P}_{r}(n_1,n_2)
&=n_1^2 \overline{P}_{r-1}(n_1,n_2)-(n_1+n_2)^2 \overline{P}_{r-1}(n_1-1,n_2), \quad r\geq 1,
\end{align*}
from which we can deduce that
$2\overline{P}_{r}(n_1,n_2)$ is divisible by $n_1{n_1+n_2\choose n_1}$ by induction on $r$.
This result is the base of the inductive proof of \cite[Theorem 1.3]{GZ10}, though it was not explicitly stated there.

\section{Proof of Theorem \ref{thm:oddmain} for $m\geq 3$}
For all nonnegative integers $a_1,\ldots,a_l$, and $k$, let
$$
C(a_1,\ldots,a_l;k)=\prod_{i=1}^l {a_i+a_{i+1}+1\choose a_i-k},
$$
where $a_{l+1}=a_1$, and let
\begin{align}\label{eq:rewriting}
S_{r}(n_1,\ldots,n_{m})
&=\frac{n_1! n_{m}!}{(n_1+n_{m}+1)!}
\sum_{k=0}^{n_1}k^r(k+1)^r(2k+1) C(n_1,\ldots,n_{m};k), \\
T_{r}(n_1,\ldots,n_{m})
&=\frac{n_1! n_{m}!}{(n_1+n_{m}+1)!}
\sum_{k=0}^{n_1}(-1)^k k^r(k+1)^r(2k+1) C(n_1,\ldots,n_{m};k).
\end{align}

Observe  that,  for $m\geq 3$, we have
\begin{align}
&\hskip -3mm C(n_1,\ldots,n_m;k)  \nonumber\\
&=\frac{(n_2+n_3+1)!(n_m+n_1+1)!}{(n_1+k+1)!(n_2-k)!(n_m+n_3+1)!}
{n_1+n_2+1\choose n_1-k}C(n_3,\ldots,n_m;k), \label{eq:cnnnk}
\end{align}
and by the Chu-Vandermonde formula (see \cite[p.~67]{AAR}) we have
\begin{align}
{n_1+n_2+1\choose n_1-k}
=\sum_{s=0}^{n_1-k}\frac{(n_1+k+1)!(n_2-k)!}{s!(s+2k+1)!(n_1-k-s)!(n_2-k-s)!}.  \label{eq:pfaff2}
\end{align}
Substituting \eqref{eq:cnnnk} and \eqref{eq:pfaff2} into the right-hand side of \eqref{eq:rewriting}, we get
\begin{align*}
S_{r}(n_1,\ldots,n_{m})&=\frac{(n_2+n_3+1)!n_1!n_m!}{(n_m+n_3+1)!}\sum_{k=0}^{n_1}\sum_{s=0}^{n_1-k}
\frac{(2k+1)^{2r+1} C(n_3,\ldots,n_{m};k)}{s!(s+2k+1)!(n_1-k-s)!(n_2-k-s)!}  \\
&=\frac{(n_2+n_3+1)!n_1!n_m!}{(n_m+n_3+1)!}\sum_{l=0}^{n_1}\sum_{k=0}^{l}
\frac{(2k+1)^{2r+1}C(n_3,\ldots,n_{m};k)}{(l-k)!(l+k+1)!(n_1-l)!(n_2-l)!},
\end{align*}
where $l=s+k$. Now, in the last sum making the substitution
$$
\frac{C(n_3,\ldots, n_{m};k)}{(l-k)!(l+k+1)!}
=\frac{(n_{m}+n_3+1)!}{(n_3+l+1)!(n_{m}+l+1)!}C(l,n_3,\ldots, n_{m};k),
$$
we obtain the following recurrence relation
\begin{align}
S_{r}(n_1,\ldots,n_{m})
=\sum_{l=0}^{n_1}{n_1\choose l}{n_2+n_3+1\choose n_2-l} S_{r}(l,n_3,\ldots,n_{m}),\quad m\geq 3. \label{eq:recsr-1}
\end{align}
Similarly, we have
\begin{align}
T_{r}(n_1,\ldots,n_{m})
=\sum_{l=0}^{n_1}{n_1\choose l}{n_2+n_3+1\choose n_2-l} T_{r}(l,n_3,\ldots,n_{m}),\quad m\geq 3. \label{eq:recsr-2}
\end{align}
We now proceed by induction on $m$. By Theorem \ref{thm:oddtwo},
we suppose that Theorem \ref{thm:oddmain} is true for $m-1$ ($m\geq 3$).
If $r=0$, then by \eqref{eq:recsr-1} and \eqref{eq:recsr-2}, Theorem \ref{thm:oddmain} is true for $m$. If $r\geq 1$,
then by definition,
$$
S_{r}(0,n_3,\ldots,n_{m})=T_{r}(0,n_3,\ldots,n_{m})=0,
$$
and by the induction hypothesis,
\begin{align*}
S_{r}(l,n_3,\ldots,n_{m}) &\equiv 0 \mod \ell^{\min\{1,r\}}n_m^{\min\{1,{r\choose 2}\}},  \\
T_{r}(l,n_3,\ldots,n_{m}) &\equiv 0 \mod \ell^{\min\{1,r\}}n_m^{\min\{1,r\}}
\end{align*}
for $\ell \geq 1$. Hence, by \eqref{eq:recsr-1} and \eqref{eq:recsr-2} and noticing that ${n_1\choose \ell}\ell=n_1{n_1-1\choose \ell-1}$, we get
\begin{align*}
S_{r}(n_1,\ldots,n_{m}) &\equiv 0 \mod n_1^{\min\{1,r\}}n_m^{\min\{1,{r\choose 2}\}},  \\
T_{r}(n_1,\ldots,n_{m}) &\equiv 0 \mod n_1^{\min\{1,r\}}n_m^{\min\{1,r\}},
\end{align*}
Namely, Theorem \ref{thm:oddmain} is true for $m\geq 3$. This completes the proof. \qed

By repeatedly using \eqref{eq:recsr-1} and \eqref{eq:recsr-2}  for $r=0,1$, we obtain the following results.
\begin{cor}\label{cor:lambda-1}
For all $m\geq 3$ and $n_1,\ldots,n_m\in\Z$, there hold
\begin{align*}
&\hskip -3mm \sum_{k=0}^{n_1}(2k+1)\prod_{i=1}^{m} {n_i+n_{i+1}+1\choose n_i-k} \\
&=(n_1+n_m+1){n_1+n_m\choose n_1}
\sum_{\lambda}{\lambda_{m-2}+n_m\choose \lambda_{m-2}}\prod_{i=1}^{m-2}{\lambda_{i-1}\choose \lambda_{i}}
{n_{i+1}+n_{i+2}+1\choose n_{i+1}-\lambda_i},  \\
&\hskip -3mm \sum_{k=0}^{n_1}k(k+1)(2k+1)\prod_{i=1}^{m} {n_i+n_{i+1}+1\choose n_i-k}  \\
&=n_m (n_1+n_m+1){n_1+n_m\choose n_1}
\sum_{\lambda}{\lambda_{m-2}+n_m-1 \choose \lambda_{m-2}-1}\prod_{i=1}^{m-2}{\lambda_{i-1}\choose \lambda_{i}}
{n_{i+1}+n_{i+2}+1\choose n_{i+1}-\lambda_i},
\end{align*}
where $n_{m+1}=\lambda_0=n_1$ and the sums are over all  sequences
$\lambda=(\lambda_1,\ldots,\lambda_{m-2})$ of nonnegative integers
such that $n_1\geq \lambda_1\geq \cdots \geq \lambda_{m-2}$.
\end{cor}

\begin{cor}\label{cor:lambda-2}
For all $m\geq 3$ and $n_1,\ldots,n_m\in\Z$, there hold
\begin{align*}
&\hskip -3mm \sum_{k=0}^{n_1}(-1)^k(2k+1)\prod_{i=1}^{m} {n_i+n_{i+1}+1\choose n_i-k} \\
&=(n_1+n_m+1){n_1+n_m\choose n_1}
\sum_{\lambda}\prod_{i=1}^{m-2}{\lambda_{i-1}\choose \lambda_{i}}
{n_{i+1}+n_{i+2}+1\choose n_{i+1}-\lambda_i},  \\
&\hskip -3mm \sum_{k=0}^{n_1}(-1)^{k+1}k(k+1)(2k+1)\prod_{i=1}^{m} {n_i+n_{i+1}+1\choose n_i-k}  \\
&={}n_m (n_1+n_m+1){n_1+n_m\choose n_1}
\sum_{\lambda}\lambda_{m-2}\prod_{i=1}^{m-2}{\lambda_{i-1}\choose \lambda_{i}}
{n_{i+1}+n_{i+2}+1\choose n_{i+1}-\lambda_i},
\end{align*}
where $n_{m+1}=\lambda_0=n_1$ and the sums are over all  sequences
$\lambda=(\lambda_1,\ldots,\lambda_{m-2})$ of nonnegative integers
such that $n_1\geq \lambda_1\geq \cdots \geq \lambda_{m-2}$.
\end{cor}

Note that the above identities show that the alternating sums in the left-hand sides are positive.
When $m=3$,  applying the Chu-Vandermonde formula, we derive the
following identities from Corollary~\ref{cor:lambda-2}:
\begin{align*}
&\hskip -3mm\sum_{k=0}^{n_1}(-1)^k (2k+1) {n_1+n_2+1\choose n_1-k}{n_2+n_3+1\choose n_2-k}{n_3+n_1+1\choose n_3-k}
=\frac{(n_1+n_2+n_3+1)!}{n_1!n_2!n_3!}, \\
&\hskip -3mm\sum_{k=0}^{n_1}(-1)^{k+1} k(k+1)(2k+1)
 {n_1+n_2+1\choose n_1-k}{n_2+n_3+1\choose n_2-k}{n_3+n_1+1\choose n_3-k}
={}\frac{(n_1+n_2+n_3)!}{(n_1-1)!(n_2-1)!(n_3-1)!}.
\end{align*}

\section{Proof of Theorem \ref{thm:oddpower}}
In this section we shall use the following theorem of Lucas (see, for example, \cite{Granville}).
We refer the reader to \cite{GZ10,SZ,ZPS} for recent applications
of Lucas' theorem.
\begin{lem}[Lucas' theorem]Let $p$ be a prime, and let $a_0,b_0,\ldots, a_m,b_m\in\{0,\ldots,p-1\}$. Then
$$
{a_0+a_1 p+\cdots + a_m p^m\choose b_0+b_1 p+\cdots + b_m p^m}\equiv \prod_{i=0}^m {a_i\choose b_i}\pmod p.
$$
\end{lem}

Suppose that $r+s\equiv 1\pmod 2$ and $s\geq 1$. Putting $m=s$ and $n_1=\cdots=n_s=n$ in Theorem \ref{thm:oddmain2}, we see that
$$
\sum_{k=0}^n (2k+1)^{r+s}{2n+1\choose n-k}^s\equiv
\sum_{k=0}^n (-1)^k(2k+1)^{r+s}{2n+1\choose n-k}^s \equiv 0 \bmod (2n+1){2n\choose n}.
$$
Therefore, by the definition of $A_{n,k}$, we have
\begin{align*}
\sum_{k=0}^{n}(2k+1)^{r}A_{n,k}^{s} \equiv \sum_{k=0}^{n}(2k+1)^{r}(-1)^k A_{n,k}^{s}
\equiv 0 \bmod \frac{{2n\choose n}}{\gcd({2n\choose n},(2n+1)^{s-1})},
\end{align*}
and so the congruence 
\eqref{eq:ank-2} holds if $s=1$.

Now suppose that $2n+1=p^a$ ($p\geq 3$) is a prime power.
By Lucas' theorem, we have
\begin{align*}
{2n\choose n}={p^a-1\choose (p^a-1)/2}\equiv {p-1\choose (p-1)/2}^a\equiv (-1)^{a(p-1)/2}\pmod p,
\end{align*}
which means that
\begin{align}
\gcd\left({2n\choose n},(2n+1)^s\right)=1. \label{eq:2n-1odd}
\end{align}
This completes the proof of Theorem \ref{thm:oddpower}.

\medskip
\noindent{\it Remark.} It is natural to wonder which numbers $n$ satisfy \eqref{eq:2n-1odd} besides those we just mentioned.
Via {\sc Maple}, we find that all such numbers $n$
less than $300$ are as follows:
\begin{align*}
&10, 27, 28, 32,37, 57, 59, 66, 85, 91, 101 ,108, 109, 118, 126, 132, 137, 150, 152, 159, \\
&164, 170, 177, 182, 188, 201, 240, 244, 252, 253, 257, 258, 271, 274, 282, 291.
\end{align*}
That is to say, Theorem \ref{thm:oddpower} is also true for these numbers $n$.
On the other hand, it is easy to see from Lucas' theorem that, if $p$ and $p^2-p+1$ are both odd primes, then
$n=(p^3-p^2+p-1)/2$ satisfies \eqref{eq:2n-1odd}. Via {\sc Maple}, we find that there are $5912$ such primes $p$ among
the first $10^5$ primes. Here we list the first $20$ such primes:
$$
3, 7, 13, 67, 79, 139, 151, 163, 193, 337, 349, 379, 457, 541, 613, 643, 727, 769, 919, 991.
$$

\section{Consequences of Theorem \ref{thm:oddmain2} and open problems}
For convenience, let $\varepsilon=\pm 1$ throughout this section. We shall give several interesting consequences of Theorem \ref{thm:oddmain2}
in this section. Note that we can also give the corresponding consequences of Theorem \ref{thm:oddmain} in the same way.

Letting $n_1=\cdots=n_{m}=n$ in Theorem~\ref{thm:oddmain2}, we have
\begin{cor}
For all $m,n\in\Z$ and $r\in\N$, there holds
$$
\sum_{k=0}^n \varepsilon^k (2k+1)^{2r+1} {2n+1\choose n-k}^{m}  \equiv 0 \mod (2n+1){2n\choose n}.
$$
\end{cor}
For $m=2a+1\geq 3$, we propose the following conjecture:
\begin{conj}
For all  $a,n\in\Z$ and $r\in\N$, there hold
\begin{align*}
{2n\choose n}^{-1}\sum_{k=0}^{n}(2k+1)^{2r+1}{2n+1\choose n-k}^{2a+1}
&\equiv
\begin{cases}
\displaystyle 1,&\text{if $n=2^b-1$}\\[5pt]
0,&\text{otherwise}
\end{cases}
\mod 2,  \\
{2n\choose n}^{-1}\sum_{k=0}^{n}(-1)^k(2k+1)^{2r+1}{2n+1\choose n-k}^{2a+1}
&\equiv
\begin{cases}
\displaystyle 1,&\text{if $n=2^b+2^c\ (b,c\in\N)$}\\[5pt]
0,&\text{otherwise}
\end{cases}
\mod 2.  \\
\end{align*}
\end{conj}
Letting $n_{2i-1}=m$ and $n_{2i}=n$ for $1\leq i\leq a$ in Theorem~\ref{thm:oddmain2},
we obtain
\begin{cor}
For all $a,m,n\in\Z$ and $r\in\N$, there holds
$$
\sum_{k=0}^m \varepsilon^k (2k+1)^{2r+1} {m+n+1\choose m-k}^a {m+n+1\choose n-k}^a  \equiv 0 \mod (m+n+1){m+n\choose m}.
$$
\end{cor}

Letting $n_{3i-2}=l$, $n_{3i-1}=m$ and $n_{3i}=n$ for $1\leq i\leq a$
in Theorem~\ref{thm:oddmain2}, we obtain
\begin{cor}
For all $a,l,m,n\in\Z$ and $r\in\N$, there holds
$$
\sum_{k=0}^l \varepsilon^k (2k+1)^{2r+1}{l+m+1\choose l-k}^a{m+n+1\choose m-k}^a {n+l+1\choose n-k}^a  \equiv 0 \mod (l+m+1){l+m\choose l}.
$$
\end{cor}

Letting $m=2a+b$, $n_1=n_3=\cdots=n_{2a-1}=n$ and letting all the other $n_i$ be $n-1$ in
Theorem~\ref{thm:oddmain2}, we get
\begin{cor}\label{cor:nn+1}For all $a,n\in\Z$ and $b,r\in \N$, there holds
\begin{align*}
\sum_{k=0}^n \varepsilon^k (2k+1)^{2r+1}{2n\choose n-k}^a{2n\choose n-k-1}^a{2n-1\choose n-k-1}^b
\equiv 0 \mod n{2n\choose n}
\end{align*}
\end{cor}

It is easy to see that Theorem \ref{thm:oddmain2} can be restated in the following form.

\begin{thm}\label{thm:rennnhalf-odd}
For all $n_1,\ldots,n_m\in\Z$ and $r\in\N$, the expression
$$
n_1!\prod_{i=1}^m\frac{(n_i+n_{i+1}+1)!}{(2n_i+1)!}
\sum_{k=0}^{n_1} \varepsilon^k (2k+1)^{2r+1}\prod_{i=1}^m {2n_i+1\choose n_i-k}\quad(n_{m+1}=-1)
$$
is always an integer.
\end{thm}

It is easy to see (cf. \cite{WZ}) that, for all $m,n\in\N$, the numbers
$\frac{(2m+1)!(2n+1)!}{(m+n+1)!m!n!}$ and $\frac{(2m)!(2n)!}{(m+n)!m!n!}$ are integers by considering
the $p$-adic order of a factorial. Letting $n_1=\cdots= n_a=m$ and
$n_{a+1}=\cdots=n_{a+b}=n$ in Theorem~\ref{thm:rennnhalf-odd}, we obtain
\begin{cor}\label{cor:mn}
For all $a,b,m,n\in\Z$ and $r\in \N$, there holds
$$
\sum_{k=0}^m \varepsilon^k (2k+1)^{2r+1}{2m+1\choose m-k}^a{2n+1\choose n-k}^b
\equiv 0 \mod \frac{(2m+1)!(2n+1)!}{(m+n+1)!m!n!}.
$$
\end{cor}

For example, we have
\begin{align}
\sum_{k=0}^n \varepsilon^k (2k+1)^{2r+1}{2n+3\choose n-k+1}^a{2n+1\choose n-k}^b
&\equiv 0 \mod (2n+3){2n+1\choose n}, \label{eq:n1+n2}\\
\sum_{k=0}^n \varepsilon^k (2k+1)^{2r+1}{4n+1\choose 2n-k}^a{2n+1\choose n-k}^b
&\equiv 0 \mod (2n+1){4n+1\choose n}, \nonumber\\
\sum_{k=0}^n \varepsilon^k (2k+1)^{2r+1}{6n+1\choose 3n-k}^a{2n+1\choose n-k}^b
&\equiv 0 \mod \frac{(6n+1)!(2n+1)!}{(4n+1)!(3n)!n!}. \nonumber
\end{align}

Similarly to the proof of Theorem \ref{thm:oddpower}, we can deduce the following results.
\begin{cor}\label{cor:n2n}
Let $r\in\N$ and $s,t\in\Z$ such that $r+s+t\equiv 1\pmod 2$, and let $n\in\Z$.
If $s=t=1$ or $\gcd\left(\frac{1}{2n+1}{2n+1\choose n},(2n+1)(2n+3)\right)=1$, then
\begin{align}
\sum_{k=0}^{n-1}\varepsilon^k (2k+1)^{r}A_{n+1,k}^s A_{n,k}^t  \equiv 0 \mod \frac{1}{2n+1}{2n+1\choose n}.
\label{eq:cor-n-1k}
\end{align}
In particular, if $2n+1$ and $2n+3$ are both prime powers, then \eqref{eq:cor-n-1k} holds.
\end{cor}
\begin{cor}
Let $r\in\N$ and $s,t\in\Z$ such that $r+s+t\equiv 1\pmod 2$, and let $n\in\Z$.
If $s=t=1$ or $\gcd\left(\frac{1}{4n+1}{4n+1\choose n},(2n+1)(4n+1)\right)=1$, then
\begin{align}
\sum_{k=0}^{n}\varepsilon^k (2k+1)^{r}A_{2n,k}^s A_{n,k}^t  \equiv 0 \mod \frac{1}{4n+1}{4n+1\choose n}.
\label{eq:cor-2nk-nk}
\end{align}
In particular, if $2n+1$ and $4n+1$ are both prime powers, then \eqref{eq:cor-2nk-nk} holds.
\end{cor}

\begin{conj}
The congruences  \eqref{eq:cor-n-1k} and \eqref{eq:cor-2nk-nk} hold for all $n\in\Z$ and $r,s,t$ given in Corollary {\rm\ref{cor:n2n}}.
\end{conj}

\begin{conj}\label{eq:twin}
There are infinitely many numbers $n\in\Z$ such that
\begin{align}
\gcd\left(\frac{1}{2n+1}{2n+1\choose n},(2n+1)(2n+3)\right)=1.  \label{eq:infinit-n-1}
\end{align}
\end{conj}

\begin{conj}There are infinitely many numbers $n\in\Z$ such that
\begin{align}
\gcd\left(\frac{1}{4n+1}{4n+1\choose n},(2n+1)(4n+1)\right)=1.  \label{eq:infinit-n-2}
\end{align}
\end{conj}

Let $f(x)$ and $g(x)$ denote the numbers of positive integers $n\leq x$
satisfying \eqref{eq:infinit-n-1} and \eqref{eq:infinit-n-2} respectively. Via {\sc Maple}, we find
that $f(x)$ and $g(x)$ grow a little slowly with respect to $x$. Table \ref{tab1} gives some values of
$f(x)$ and $g(x)$.

\begin{table}[h]
\begin{center}{\footnotesize
\begin{tabular}{|c|c|c|c|c|c|c|c|c|c|c|c|c|} \hline
$n$    & 1 & 10 &20 &50 & 100 &200 & 500 & 1000 & 2000 &3000& 4000 &5000 \\ \hline
$f(n)$ & 1 & 8  &13 &24 & 38  &59  & 104 & 167  & 255  &353 & 439  &508  \\ \hline
$g(n)$ & 1 & 7  &14 &23 & 37  &56  & 108 & 169  & 270  &366 & 445  &523  \\ \hline
\end{tabular} }
\end{center}
\caption{Some values of $f(x)$ and $g(x)$.\label{tab1}}
\end{table}

Using Lucas' theorem, it is easy to see that if $2n+1$ and $2n+3$ are both prime powers, then \eqref{eq:infinit-n-1} holds.
It is well known that the {\it twin prime conjecture} states that there are infinitely many primes $p$ such that $p+2$ is also prime.
Therefore, if the twin prime conjecture is true, then so is Conjecture \ref{eq:twin}. Since the twin prime conjecture
is rather difficult, we hope that Conjecture \ref{eq:twin} might be tackled in another way.

In fact, we have the following more general conjecture:
\begin{conj}For all $r\in\N$ and $m,n,s,t\in\Z$ such that $r+s+t\equiv 1\pmod 2$, there holds
$$
(m+n+1)\sum_{k=0}^m \varepsilon^k (2k+1)^r A_{m,k}^s A_{n,k}^t
\equiv 0 \mod \frac{(2m)!(2n)!}{(m+n)!m!n!}.
$$
\end{conj}
Note that the numbers $\frac{(2m)!(2n)!}{(m+n)!m!n!}$ are called super Catalan numbers,
of which no combinatorial interpretations are known for general $m$ and $n$ until now.
See Gessel \cite{Gessel} and Georgiadis et al. \cite{GMT}.

{}From Theorem~\ref{thm:rennnhalf-odd} it is easy to see that, for all $a_1,\ldots,a_m\in\Z$,
\begin{align}
n_1!\prod_{i=1}^m\frac{(n_i+n_{i+1}+1)!}{(2n_i+1)!}
\sum_{k=0}^{n_1} \varepsilon^k (2k+1)^{2r+1}\prod_{i=1}^m {2n_i+1\choose n_i-k}^{a_i}\quad(n_{m+1}=-1) \label{eq:final-form}
\end{align}
is an integer.
For $m=3$, letting $(n_1,n_2,n_3)$ be $(n,n+2,n+1)$, $(n,3n,2n)$, $(2n,n,3n)$,
$(2n,n,4n)$, or $(3n,2n,4n)$, we obtain the following three corollaries.
\begin{cor}\label{cor:n1n2n3}
For all $a,b,c,n\in\Z$ and $r\in \N$, there holds
\begin{align}
\sum_{k=0}^n \varepsilon^k (2k+1)^{2r+1}   {2n+1\choose n-k}^a {2n+3\choose n-k+1}^b {2n+5\choose n-k+2}^c
&\equiv 0 \mod (2n+5){2n+1\choose n}. \label{eq:n1n2n3}
\end{align}
\end{cor}
\begin{cor}\label{cor:rst-246n}
For all $a,b,c,n\in\Z$ and $r\in \N$, there hold
\begin{align*}
\sum_{k=0}^n \varepsilon^k (2k+1)^{2r+1} {6n+1\choose 3n-k}^a {4n+1\choose 2n-k}^b {2n+1\choose n-k}^c
&\equiv 0 \mod (2n+1){6n+1\choose n}, \\
\sum_{k=0}^n \varepsilon^k (2k+1)^{2r+1} {6n+1\choose 3n-k}^a {4n+1\choose 2n-k}^b {2n+1\choose n-k}^c
&\equiv 0 \mod (2n+1){6n+1\choose 3n}.
\end{align*}
\end{cor}

\begin{cor}\label{cor:rst-248n}
For all $a,b,c,n\in\Z$ and $r\in \N$, there hold
\begin{align*}
2\sum_{k=0}^n \varepsilon^k (2k+1)^{2r+1} {8n+1\choose 4n-k}^a {4n+1\choose 2n-k}^b {2n+1\choose n-k}^c
&\equiv 0 \mod \frac{2(2n+1)(4n+1)}{3n+1}{8n+1\choose 3n},   \\
\sum_{k=0}^n  \varepsilon^k (2k+1)^{2r+1} {8n+1\choose 4n-k}^a {6n+1\choose 3n-k}^b {4n+1\choose 2n-k}^c
&\equiv 0 \mod (4n+1){8n+1\choose 3n},
\end{align*}
\end{cor}
Note that
$$
\frac{2}{3n+1}{8n+1\choose 3n}=8{8n+1\choose 3n}-3{8n+2\choose 3n+1}\in\Z.
$$

\begin{conj}
For all $n,r,s,t\in\Z$ such that $r+s+t\equiv 1\pmod 2$, there hold
\begin{align*}
(4n+1)\sum_{k=0}^n \varepsilon^k A_{3n,k}^r A_{2n,k}^s A_{n,k}^t
&\equiv 0 \mod \frac{1}{6n+1}{6n+1\choose n}, \\
(4n+1)\sum_{k=0}^n \varepsilon^k A_{3n,k}^r A_{2n,k}^s A_{n,k}^t
&\equiv 0 \mod \frac{1}{6n+1}{6n+1\choose 3n}, \\
\sum_{k=0}^n \varepsilon^k A_{4n,k}^r A_{2n,k}^s A_{n,k}^t
&\equiv 0 \mod \frac{1}{8n+1}{8n+1\choose 3n}, \\
(6n+1)\sum_{k=0}^n \varepsilon^k A_{4n,k}^r A_{3n,k}^s A_{2n,k}^t
&\equiv 0 \mod \frac{1}{8n+1}{8n+1\choose 3n}.
\end{align*}
\end{conj}

Finally, for general $m\geq 2$, in \eqref{eq:final-form} taking $(n_1,\ldots,n_{m})$ to be
\begin{align*}
\begin{cases}
(n,n+2,n+4,\ldots,n+m-1,n+m-2,n+m-4,n+m-6,\ldots,n+1),&\text{if $m$ is odd,}\\[5pt]
(n+1,n+3,n+5,\ldots,n+m-1,n+m-2,n+m-4,n+m-6,\ldots,n),&\text{if $m$ is even,}
\end{cases}
\end{align*}
we get the following generalization of congruences \eqref{eq:n1+n2} and \eqref{eq:n1n2n3}.
\begin{cor}\label{cor:final}
Let $m\geq 2$, and let $n,a_1,\ldots,a_m\in\Z$ and $r\in\N$. Then
\begin{align*}
\sum_{k=0}^n \varepsilon^k (2k+1)^{2r+1} \prod_{i=1}^m {2n+2i-1\choose n+i-k-1}^{a_i}
&\equiv 0 \mod (2n+2m-1){2n+1\choose n}.
\end{align*}
\end{cor}
We end this paper with the following challenging conjecture related to Corollary \ref{cor:final}.
\begin{conj}\label{conj:final}
For all $n,r_1,\ldots,r_m\in\Z$ such that $r_1+\cdots+r_m\equiv 1\pmod 2$, there holds
\begin{align*}
\sum_{k=0}^n \varepsilon^k \prod_{i=1}^m A_{n+i-1,k}^{r_i}
&\equiv 0 \mod \frac{1}{2n+1}{2n+1\choose n}.
\end{align*}
\end{conj}
Note that, the $m=1$ case of Conjecture \ref{conj:final} is a little
weaker than our previous Conjecture \ref{conj:first}, where the
modulus $\frac{1}{2n+1}{2n+1\choose n}$ is replaced by ${2n\choose
n}$. By Corollary \ref{cor:final}, it is easy to see that Conjecture
\ref{conj:final} is true for $n=2$, or $m\leq 6$ and
$n=4,9,10,11,3280,7651,7652$ (with the help of {\sc Maple}).

\vskip 5mm
\noindent{\bf Acknowledgments.} This work was partially
supported by the Fundamental Research Funds for the Central Universities,  Shanghai Rising-Star Program (\#09QA1401700),
Shanghai Leading Academic Discipline Project (\#B407), and the National Science Foundation of China (\#10801054).

\renewcommand{\baselinestretch}{1}

\end{document}